\documentclass[a4paper,draft]{amsproc}
\usepackage{amsfonts}
\usepackage{amssymb}
\usepackage{amscd}

\theoremstyle{plain}

\theoremstyle{definition}

\theoremstyle{remark}

 \numberwithin{equation}{section}
\renewcommand{\leq}{\leqslant}

\email{gulnur.saffak@omu.edu.tr, erginbayram@yahoo.com, kasape@omu.edu.tr}
\input{tcilatex}

\begin{document}
\title[Surfaces with a common asymptotic in $R_{1}^{3}$]{SURFACES WITH A
COMMON ASYMPTOTIC CURVE IN MINKOWSKI 3-SPACE}
\author[\c{S}affak, Bayram and Kasap]{\bfseries G\"{u}lnur \c{S}affak, Ergin
Bayram and Emin Kasap}
\address{ Department of Mathematics \\
Faculty of Arts and Sciences\\
Ondokuz Mayis University \\
Samsun\\
Turkey}

\begin{abstract}
In this paper, we express surfaces parametrically through a given spacelike
(timelike) asymptotic curve using the Frenet frame of the curve in Minkowski
3-space. Necessary and sufficient conditions for the coefficients of the
Frenet frame to satisfy both parametric and asymptotic requirements are
derived. We also present some interesting examples to show the validity of
this study.
\end{abstract}

\maketitle

\vspace{18mm} \setcounter{page}{1} \thispagestyle{empty}

\section{Introduction}

Tangent
vectors on a Minkowski surface are classified into three types. A smooth
curve on a surface is said to be timelike, null, or spacelike if its tangent
vectors are always timelike, null or spacelike, respectively. Physically, a
timelike curve corresponds to the path of an observer moving at less than
the speed of light. Null curves correspond to moving at the speed of light,
and spacelike curves to moving faster than light.

The concept of family of surfaces having a given characteristic curve was
first introduced by Wang et.al. \cite{wang} in Euclidean 3-space. Kasap
et.al. \cite{kasap akyildiz orbay} generalized the work of Wang by
introducing new types of marching-scale functions, coefficients of the
Frenet frame appearing in the parametric representation of surfaces. In \cite%
{kasap akyildiz} Kasap and Aky\i ld\i z defined surfaces with a common
geodesic in Minkowski 3-space and gave the sufficient conditions on
marching-scale functions so that the given curve is a common geodesic on
that surfaces. \c{S}affak and Kasap \cite{saffak kasap} studied family of
surfaces with a common null geodesic.

With the inspiration of work of Wang, Li et.al.\cite{Li} changed the
characteristic curve from geodesic to line of curvature and defined the
surface pencil with a common line of curvature. Recently, in \cite{bayram}
Bayram et.al. defined the surface pencil with a common asymptotic curve.
They introduced three types of marching-scale functions and derived the
necessary and sufficient conditions on them to satisfy both parametric and
asymptotic requirements.

In this paper, we study the problem: given a 3D spacelike (timelike) curve,
how to characterize those surfaces that posess this curve as a common
parametric and asymptotic curve in Minkowski 3-space. In section 2, we give
some preliminary information about curves and surfaces in Minkowski 3-space
and define isoasymptotic curve. We express spacelike surfaces as a linear
combination of the Frenet frame of the given curve and derive necessary and
sufficient conditions on marching-scale functions to satisfy both parametric
and asymptotic requirements in Section 3. Section 4 is devoted to timelike
surfaces. We illustrate the method by giving some examples. Also, all
minimal timelike surfaces are given as examples of timelike surfaces with
common asymptotic curve.

\section{Preliminaries}

Let us consider Minkowski 3-space $%
\mathbb{R}
_{1}^{3}=\left[ 
\mathbb{R}
^{3},\left( +,+,-\right) \right] $ and let the Lorentzian inner product of $%
X=\left( x_{1},x_{2},x_{3}\right) $ and $Y=\left( y_{1},y_{2},y_{3}\right)
\in 
\mathbb{R}
_{1}^{3}$ be 
\begin{equation*}
\left\langle X,Y\right\rangle =x_{1}y_{1}+x_{2}y_{2}-x_{3}y_{3.}
\end{equation*}

A vector $X\in 
\mathbb{R}
_{1}^{3}$ is called a spacelike vector when $\left\langle X,X\right\rangle
>0 $ or $X=0.$ It is called timelike and null (light-like) vector in case of 
$\left\langle X,X\right\rangle <0,$ and $\left\langle X,X\right\rangle =0~$%
for $X\neq 0,~$respectively, \cite{oneill}.

The vector product of vectors $X=\left( x_{1},x_{2},x_{3}\right) \ $and $%
Y=\left( y_{1},y_{2},y_{3}\right) \ $in$\ 
\mathbb{R}
_{1}^{3}~$is defined by \cite{agutagawa}

\begin{equation*}
XxY=\left(
x_{3}y_{2}-x_{2}y_{3},x_{1}y_{3}-x_{3}y_{1},x_{1}y_{2}-x_{2}y_{1}\right) .
\end{equation*}

Let $\alpha =\alpha \left( u\right) ~$be a unit speed curve in $%
\mathbb{R}
_{1}^{3}.$ We denote the natural curvature and torsion of $\alpha \left(
u\right) ~$with $\kappa \left( u\right) $\ and $\tau \left( u\right) $,
respectively. Consider the Frenet frame $\left\{ e_{1},e_{2},e_{3}\right\} $
associated with curve $\alpha \left( u\right) $ such that $e_{1}=e_{1}\left(
u\right) ,\ e_{2}=e_{2}\left( u\right) $ and $e_{3}=e_{3}\left( u\right) ~$%
are the unit tangent, the principal normal and the binormal vector fields,
respectively. If $\alpha =\alpha \left( u\right) $ is a spacelike curve,
then the structural equations (or Frenet formulas) of this frame given as%
\begin{equation*}
e_{1}^{\prime }\left( u\right) =\kappa \left( u\right) e_{2}\left( u\right)
,\ e_{2}^{\prime }\left( u\right) =\varepsilon \kappa \left( u\right)
e_{1}\left( u\right) +\tau \left( u\right) e_{3},\ e_{3}^{\prime }\left(
u\right) =\tau \left( u\right) e_{2}\left( u\right) ,
\end{equation*}%
where 
\begin{equation*}
\varepsilon =\left\{ 
\begin{array}{c}
-1,~e_{3}\ is\ timelike, \\ 
1,\ e_{3}\ is\ spacelike.%
\end{array}%
\right.
\end{equation*}

If $\alpha =\alpha \left( u\right) $ is a timelike curve, then above
equations are given as \cite{woestijne}

\begin{equation*}
e_{1}^{\prime }\left( u\right) =\kappa \left( u\right) e_{2}\left( u\right)
,\ e_{2}^{\prime }\left( u\right) =\kappa \left( u\right) e_{1}\left(
u\right) -\tau \left( u\right) e_{3},\ e_{3}^{\prime }\left( u\right) =\tau
\left( u\right) e_{2}\left( u\right) .
\end{equation*}

A surface in $%
\mathbb{R}
_{1}^{3}$ is called a timelike surface if the induced metric on the surface
is a Lorentz metric and is called a spacelike surface if the induced metric
on the surface is a positive definite Riemannian metric, i.e. the normal
vector on the spacelike (timelike) surface is a timelike (spacelike) vector 
\cite{beem}.

A curve on a surface is called an asymptotic curve $\alpha \left( u\right) $
provided its velocity always points in an asymptotic direction, that is the
direction in which the normal curvature is zero \cite{duggal}. According to
the above definition the curve is an asymptotic curve on the surface $%
\varphi \left( u,v\right) $ if and only if%
\begin{equation}
\frac{\partial n\left( u,v_{0}\right) }{\partial u}\cdot e_{1}\left(
u\right) =0  \label{2.1}
\end{equation}%
where \textquotedblleft $\cdot $\textquotedblright\ denotes the Lorentzian
inner product and n is a normal vector of $\varphi =\varphi \left(
u,v\right) $ \cite{duggal}.

An isoparametric curve $\alpha =\alpha \left( u\right) $ is a curve on a
surface $\varphi =\varphi \left( u,v\right) $ in $%
\mathbb{R}
_{1}^{3}$ is that has a constant u or v - parameter value. In other words,
there exists a parameter $u_{0}$\ or $v_{0}$ such that $\alpha \left(
u\right) =\varphi \left( u,v_{0}\right) $ or $\alpha \left( v\right)
=\varphi \left( u_{0},v\right) $. Given a parametric curve $\alpha \left(
u\right) $, we call it an \textit{isoasymptotic} of the surface $\varphi $
if it is both an asymptotic curve and a parameter curve on $\varphi $.

We assume that $\kappa \left( u\right) \neq 0$ for $\alpha \left( u\right) $
along the paper. Otherwise, the principal normal of the curve is undefined
or the curve is a straightline.

\section{Spacelike surfaces with a common spacelike asymptotic}

Let $\varphi =\varphi \left( u,v\right) $ be a parametric spacelike surface.
The surface is defined by a given curve $\alpha =\alpha \left( u\right) $ as
follows:

\begin{equation}
\left\{ 
\begin{array}{c}
\varphi \left( u,v\right) =\alpha \left( u\right) +\left[ x\left( u,v\right)
e_{1}\left( u\right) +y\left( u,v\right) e_{2}\left( u\right) +z\left(
u,v\right) e_{3}\left( u\right) \right] ,\  \\ 
L_{1}\leq u\leq L_{2},\ T_{1}\leq v\leq T_{2},%
\end{array}%
\right.  \label{3.1}
\end{equation}%
where $x\left( u,v\right) ,\ y\left( u,v\right) \ $and $z\left( u,v\right) \ 
$are $C^{1}$ functions.The values of the marching-scale functions $x\left(
u,v\right) ,\ y\left( u,v\right) \ $and $z\left( u,v\right) $ indicate,
respectively; the extension-like, flexion-like and retortion-like effects,
by the point unit through the time $v$, starting from $\alpha =\alpha \left(
u\right) $ and $\left\{ e_{1}\left( u\right) ,\ e_{2}\left( u\right) ,\
e_{3}\left( u\right) \right\} $ is the Frenet frame associated with the
curve $\alpha \left( u\right) .\ $

Our goal is to find the necessary and sufficient conditions for which the
curve is a parameter curve and an asymptotic curve on the surface .

Firstly, since $\alpha \left( u\right) $ is a parameter curve on the surface 
$\varphi \left( u,v\right) $, there exists a parameter $v_{0}\in \left[
T_{1},T_{2}\right] $ such that

\begin{equation*}
x\left( u,v_{0}\right) =y\left( u,v_{0}\right) =z\left( u,v_{0}\right) =0,\
L_{1}\leq u\leq L_{2},\ T_{1}\leq v\leq T_{2}.
\end{equation*}

Secondly, since $\alpha \left( u\right) $ is an asymptotic curve on the
surface, from Eqn. \ref{2.1} there exists a parameter $v_{0}\in \left[
T_{1},T_{2}\right] $ such that $\frac{\partial n}{\partial u}\left(
u,v_{0}\right) \cdot e_{1}\left( u\right) =0.$

\begin{theorem}
A spacelike curve $\alpha \left( u\right) $ is isoasymptotic on a spacelike
surface $\varphi \left( u,v\right) $ if and only if the following conditions
are satisfied: 
\begin{equation}
\left\{ 
\begin{array}{c}
x\left( u,v_{0}\right) =y\left( u,v_{0}\right) =z\left( u,v_{0}\right) =0,
\\ 
\frac{\partial z}{\partial v}\left( u,v_{0}\right) =0.%
\end{array}%
\right.   \label{3.3}
\end{equation}

\begin{proof}
Let $\alpha \left( u\right) $ be a spacelike curve on a spacelike surface $%
\varphi \left( u,v\right) $. If $\alpha \left( u\right) $ is a parameter
curve on this surface, then there exists a parameter $v=v_{0}$ such that $%
\alpha \left( u\right) =\varphi \left( u,v_{0}\right) $, that is,%
\begin{equation}
x\left( u,v_{0}\right) =y\left( u,v_{0}\right) =z\left( u,v_{0}\right) =0.
\label{3.4}
\end{equation}

The normal vector of $\varphi =\varphi \left( u,v\right) \ $can be written as%
\begin{equation*}
n\left( u,v\right) =\frac{\partial \varphi \left( u,v\right) }{\partial u}x%
\frac{\partial \varphi \left( u,v\right) }{\partial v}
\end{equation*}%
Since%
\begin{eqnarray*}
\frac{\partial \varphi \left( u,v\right) }{\partial u} &=&\left( 1+\frac{%
\partial x\left( u,v\right) }{\partial u}+\varepsilon \kappa \left( u\right)
y\left( u,v\right) \right) e_{1}\left( u\right)  \\
&&+\left( \frac{\partial y\left( u,v\right) }{\partial u}+\kappa \left(
u\right) x\left( u,v\right) +\tau \left( u\right) z\left( u,v\right) \right)
e_{2}\left( u\right)  \\
&&+\left( \frac{\partial z\left( u,v\right) }{\partial u}+\tau \left(
u\right) y\left( u,v\right) \right) e_{3}\left( u\right) ,
\end{eqnarray*}%
\begin{equation*}
\frac{\partial \varphi \left( u,v\right) }{\partial v}=\frac{\partial
x\left( u,v\right) }{\partial v}e_{1}\left( u\right) +\frac{\partial y\left(
u,v\right) }{\partial v}e_{2}\left( u\right) +\frac{\partial z\left(
u,v\right) }{\partial v}e_{3}\left( u\right) ,
\end{equation*}%
the normal vector can be expressed as%
\begin{eqnarray*}
n\left( u,v\right)  &=&\left[ \left( \frac{\partial y\left( u,v\right) }{%
\partial u}+\kappa \left( u\right) x\left( u,v\right) +\tau \left( u\right)
z\left( u,v\right) \right) \frac{\partial z\left( u,v\right) }{\partial v}%
\right.  \\
&&\left. -\left( \frac{\partial z\left( u,v\right) }{\partial u}+\tau \left(
u\right) y\left( u,v\right) \right) \frac{\partial y\left( u,v\right) }{%
\partial v}\right] e_{1} \\
&&+\left[ \left( 1+\frac{\partial x\left( u,v\right) }{\partial u}%
+\varepsilon \kappa \left( u\right) y\left( u,v\right) \right) \frac{%
\partial z\left( u,v\right) }{\partial v}\right.  \\
&&\left. -\left( \frac{\partial z\left( u,v\right) }{\partial u}+\tau \left(
u\right) y\left( u,v\right) \right) \frac{\partial x\left( u,v\right) }{%
\partial v}\right] e_{2} \\
&&+\left[ \left( 1+\frac{\partial x\left( u,v\right) }{\partial u}%
+\varepsilon \kappa \left( u\right) y\left( u,v\right) \right) \frac{%
\partial y\left( u,v\right) }{\partial v}\right.  \\
&&-\left. \left( \frac{\partial y\left( u,v\right) }{\partial u}+\kappa
\left( u\right) x\left( u,v\right) +\tau \left( u\right) z\left( u,v\right)
\right) \frac{\partial x\left( u,v\right) }{\partial v}\right] e_{3}
\end{eqnarray*}%
Thus, if we let%
\begin{equation*}
\phi _{1}\left( u,v_{0}\right) =\frac{\partial y\left( u,v_{0}\right) }{%
\partial u}\frac{\partial z\left( u,v_{0}\right) }{\partial v}-\frac{%
\partial z\left( u,v_{0}\right) }{\partial u}\frac{\partial y\left(
u,v_{0}\right) }{\partial v},
\end{equation*}%
\begin{equation*}
\phi _{2}\left( u,v_{0}\right) =\left( 1+\frac{\partial x\left(
u,v_{0}\right) }{\partial u}\right) \frac{\partial z\left( u,v_{0}\right) }{%
\partial v}-\frac{\partial z\left( u,v_{0}\right) }{\partial u}\frac{%
\partial x\left( u,v_{0}\right) }{\partial v},
\end{equation*}%
\begin{equation*}
\phi _{3}\left( u,v_{0}\right) =\left( 1+\frac{\partial x\left(
u,v_{0}\right) }{\partial u}\right) \frac{\partial y\left( u,v_{0}\right) }{%
\partial v}-\frac{\partial y\left( u,v_{0}\right) }{\partial u}\frac{%
\partial x\left( u,v_{0}\right) }{\partial v},
\end{equation*}%
we obtain%
\begin{equation*}
n\left( u,v_{0}\right) =\phi _{1}\left( u,v_{0}\right) e_{1}\left( u\right)
+\phi _{2}\left( u,v_{0}\right) e_{2}\left( u\right) +\phi _{3}\left(
u,v_{0}\right) e_{3}\left( u\right) .
\end{equation*}%
From Eqn. \ref{2.1}, we know that $\alpha \left( u\right) $ is an asymptotic
curve if and only if%
\begin{equation*}
\frac{\partial \phi _{1}\left( u,v_{0}\right) }{\partial u}+\kappa \left(
u\right) \phi _{2}\left( u,v_{0}\right) =0.
\end{equation*}%
Since $\kappa \left( u\right) =\left\Vert \alpha ^{\prime \prime }\left(
u\right) \right\Vert \neq 0,\ \phi _{2}\left( u,v_{0}\right) =\frac{\partial
z\left( u,v_{0}\right) }{\partial v}\ $and by Eqn. \ref{3.4} we have $\frac{%
\partial \phi _{1}\left( u,v_{0}\right) }{\partial u}=0.$ Therefore, Eqn. \ref%
{2.1} is simplified to%
\begin{equation*}
\frac{\partial z\left( u,v_{0}\right) }{\partial v}=0,
\end{equation*}%
which completes the proof.
\end{proof}
\end{theorem}

We call the set of surfaces defined by Eqn. \ref{3.1} and satisfying Eqn. %
\ref{3.3} \textit{the family of spacelike surfaces with a common spacelike
asymptotic}. Any surface $\varphi \left( u,v\right) $ defined by Eqn. \ref%
{3.1} and satisfying Eqn. \ref{3.3} is a member of this family.

In Eqn. \ref{3.1}, marching-scale functions $x\left( u,v\right) ,\ y\left(
u,v\right) $ and $z\left( u,v\right) \ $can be choosen in two different
forms:

\textbf{1)}\qquad If we choose%
\begin{equation}
\left\{ 
\begin{array}{c}
x\left( u,v\right) =\sum_{i=1}^{p}a_{1i}l\left( u\right) ^{i}X\left(
v\right) ^{i}, \\ 
y\left( u,v\right) =\sum_{i=1}^{p}a_{2i}m\left( u\right) ^{i}Y\left(
v\right) ^{i}, \\ 
z\left( u,v\right) =\sum_{i=1}^{p}a_{3i}l\left( u\right) ^{i}Z\left(
v\right) ^{i},%
\end{array}%
\right.  \label{3.5}
\end{equation}%
then we can simply express the sufficient condition for which the spacelike
curve $\alpha \left( u\right) $ is an isoasymptotic on the spacelike surface 
$\varphi \left( u,v\right) $ as%
\begin{equation}
\left\{ 
\begin{array}{c}
X\left( v_{0}\right) =Y\left( v_{0}\right) =Z\left( v_{0}\right) =0, \\ 
a_{31}=0\ or\ n\left( u\right) =0\ or\ \frac{dZ\left( v_{0}\right) }{dv}=0,%
\end{array}%
\right.  \label{3.6}
\end{equation}%
where $l\left( u\right) ,\ m\left( u\right) ,\ n\left( u\right) ,\ X\left(
v\right) ,\ Y\left( v\right) \ $and $Z\left( v\right) \ $are $C^{1}\ $%
functions and $a_{ij}\in 
\mathbb{R}
,\ i=1,\ 2,\ 3,\ \ j=1,\ 2,...,\ p.$

\textbf{2) }If we choose%
\begin{equation}
\left\{ 
\begin{array}{c}
x\left( u,v\right) =f\left( \sum_{i=1}^{p}a_{1i}l\left( u\right) ^{i}X\left(
v\right) ^{i}\right) , \\ 
y\left( u,v\right) =g\left( \sum_{i=1}^{p}a_{2i}m\left( u\right) ^{i}Y\left(
v\right) ^{i}\right) , \\ 
z\left( u,v\right) =h\left( \sum_{i=1}^{p}a_{3i}l\left( u\right) ^{i}Z\left(
v\right) ^{i}\right) ,%
\end{array}%
\right.  \label{3.7}
\end{equation}%
then we can simply express the sufficient condition for which the spacelike
curve is an isoasymptotic on the spacelike surface $\varphi \left(
u,v\right) $ as%
\begin{equation}
\left\{ 
\begin{array}{c}
X\left( v_{0}\right) =Y\left( v_{0}\right) =Z\left( v_{0}\right) =0\ and\
f\left( 0\right) =g\left( 0\right) =h\left( 0\right) =0, \\ 
a_{31}=0\ or\ n\left( u\right) =0\ or\ h^{\prime }\left( 0\right) =0\ or\ 
\frac{dZ\left( v_{0}\right) }{dv}=0,%
\end{array}%
\right.  \label{3.8}
\end{equation}%
where $l\left( u\right) ,\ m\left( u\right) ,\ n\left( u\right) ,\ X\left(
v\right) ,\ Y\left( v\right) \ $and $Z\left( v\right) ,\ \ f\ ,\ g\ $and $h$%
\ are $C^{1}\ $functions and $a_{ij}\in 
\mathbb{R}
,\ i=1,\ 2,\ 3,\ \ j=1,\ 2,...,\ p.$\bigskip

Because the parameters $a_{ij},\ i=1,\ 2,\ 3,\ \ j=1,\ 2,...,\ p$ in Eqns. %
\ref{3.5} and \ref{3.7} control the shape of the surface, one can adjust
these parameters to produce spacelike surfaces that meet certain
constraints, such as conditions on the boundary, curvature, etc. The
marching-scale functions in Eqns. \ref{3.5} and \ref{3.7} are general for
expressing surfaces with a given curve as an isoasymptotic curve.
Furthermore, conditions for different types of marching-scale functions can
be obtained from Eqn. \ref{3.3}.

Because there are no constraints related to the given curve in Eqns \ref{3.6}
or \ref{3.8}, a spacelike surface family passing through a given regular arc
length curve, acting as both a parametric curve and an asymptotic curve, can
always be found by choosing suitable marching-scale functions.

\begin{example}
\label{example 3.1}Suppose we are given a parametric spacelike curve 
\begin{equation*}
\alpha \left( u\right) =\left( \cos u,\sin u,0\right) ,0\leq u\leq 2\pi .
\end{equation*}%
We will construct a family of spacelike surfaces sharing the curve $\alpha
\left( u\right) $ as the spacelike isoasymptotic. It is easy to show that%
\begin{equation*}
\left\{ 
\begin{array}{c}
e_{1}\left( u\right) =\left( -\sin u,\cos u,0\right) , \\ 
e_{2}\left( u\right) =\left( -\cos u,-\sin u,0\right) , \\ 
e_{3}=\left( 0,0,1\right) .%
\end{array}%
\right. 
\end{equation*}%
If we choose 
\begin{equation*}
x\left( u,v\right) =0,\ y\left( u,v\right) =\cos v+\sum_{k=2}^{p}a_{2k}\cos
^{k}v,\ z\left( u,v\right) =\sum_{k=1}^{p}a_{3k}\left( 1+\sin v\right) ^{k}\ 
\end{equation*}%
and $v_{0}=\frac{3\pi }{2}$ then Eqn.\ref{3.6} is satisfied. Thus, we
immediately obtain a member of this family as%
\begin{eqnarray*}
\varphi \left( u,v\right)  &=&\left( \cos u\left( 1-\cos v-\sum_{k=2}^{4}%
\frac{1}{2}\left( \cos v\right) ^{k}\right) ,\right.  \\
&&\sin u\left( 1-\cos v-\sum_{k=2}^{4}\frac{1}{2}\left( \cos v\right)
^{k}\right) , \\
&&\left. \sum_{k=1}^{4}\frac{1}{2}\left( 1+\sin v\right) ^{k}\right) ,
\end{eqnarray*}%
where $4\leq v\leq 5$ (Figure 1).
\begin{figure}[tbh]
\caption{{}A member of spacelike surface family and its common spacelike
asymptotic.}
\end{figure}
\end{example}

\section{Timelike surfaces with a common spacelike or timelike asymptotic}

Let $\varphi \left( u,v\right) $ be a parametric timelike surface. The
surface is defined by a given curve $\alpha =\alpha \left( u\right) $ as
follows:%
\begin{equation}
\left\{ 
\begin{array}{c}
\varphi \left( u,v\right) =\alpha \left( u\right) +\left[ x\left( u,v\right)
e_{1}\left( u\right) +y\left( u,v\right) e_{2}\left( u\right) +z\left(
u,v\right) e_{3}\left( u\right) \right] ,\  \\ 
L_{1}\leq u\leq L_{2},\ T_{1}\leq v\leq T_{2},%
\end{array}%
\right.  \label{4.1}
\end{equation}%
where $x\left( u,v\right) ,\ y\left( u,v\right) \ $and $z\left( u,v\right) \ 
$are $C^{1}$ functions and $\left\{ e_{1}\left( u\right) ,e_{2}\left(
u\right) ,e_{3}\left( u\right) \right\} \ $is the Frenet frame associated
with the curve $\alpha \left( u\right) .$

Similar computation shows that the conditions \ref{3.3}, \ref{3.6} and \ref%
{3.8} are valid for a curve to be both isoparametric and asymptotic on
timelike surfaces. We call the set of surfaces defined by Eqn. \ref{4.1} and
satisfying Eqn. \ref{3.3} \textit{the family of timelike surfaces with a
common timelike asymptotic}. Any surface defined by Eqn. \ref{4.1} and
satisfying Eqn. \ref{3.3} is a member of this family.

Now let us give some examples for timelike surfaces with a common asymptotic
curve (spacelike or timelike):

\begin{example}
\label{example 4.1}Suppose we are given a parametric timelike curve 
\begin{equation*}
\alpha \left( u\right) =\left( \cosh u,0,\sinh u\right), 
\end{equation*}
where $-2\leq u\leq 2$. We will construct a family of timelike surfaces
sharing the curve $\alpha \left( u\right) $ as the timelike isoasymptotic.
It is easy to show that%
\begin{equation*}
\left\{ 
\begin{array}{c}
e_{1}\left( u\right) =\left( \sinh u,0,\cosh u\right) , \\ 
e_{2}\left( u\right) =\left( \cosh u,0,\sinh u\right) , \\ 
e_{3}=\left( 0,1,0\right) .%
\end{array}%
\right. 
\end{equation*}%
If we choose%
\begin{equation*}
x\left( u,v\right) =0,\ y\left( u,v\right) =\sin v,\ z\left( u,v\right)
=uv^{2}\ 
\end{equation*}%
and $v_{0}=0$ then Eqn.\ref{3.6} is satisfied. By putting these functions
into Eqn. \ref{4.1}, we obtain the following timelike surface passing
through the common asymptotic $\alpha \left( u\right) $%
\begin{equation*}
\varphi \left( u,v\right) =\left( \cosh u+\sin v\cosh u,uv^{2},\sinh u+\sin
v\sinh u\right), 
\end{equation*}%
where $-1\leq v\leq 1\ \left( Figure\ 2\right) .$%
\begin{figure}[tbh]
\caption{{}A member of timelike surface family and its common timelike
asymptotic.}
\end{figure}
\end{example}

\begin{example}
\label{example 4.2}Suppose we are given a parametric spacelike curve 
\begin{equation*}
\alpha \left( u\right) =\left( \cos u,\sin u,0\right), 
\end{equation*}
where 0$\leq u\leq 2\pi $. We will consruct a family of spacelike surfaces
sharing the curve $\alpha \left( u\right) $ as a spacelike isoasymptotic. It
is easy to show that%
\begin{equation*}
\left\{ 
\begin{array}{c}
e_{1}\left( u\right) =\left( -\sin u,\cos u,0\right) , \\ 
e_{2}\left( u\right) =\left( -\cos u,-\sin u,0\right) , \\ 
e_{3}=\left( 0,0,1\right) .%
\end{array}%
\right. 
\end{equation*}%
By choosing marching-scale functions as%
\begin{eqnarray*}
x\left( u,v\right)  &=&0,\ y\left( u,v\right) =\sin \left(
\sum_{k=1}^{p}a_{2k}\left( \sinh v\right) ^{k}\right) ,\  \\
z\left( u,v\right)  &=&\sin \left( \sum_{k=1}^{p}a_{3k}\left( 1-\cosh
v\right) ^{k}\right) ,
\end{eqnarray*}%
where $a_{2k},a_{3k}\in 
\mathbb{R}
\ $and letting $v_{0}=0,$ then Eqn. \ref{3.8} is satisfied. Thus, we obtain a
member of the surface family with a common spacelike asymptotic curve $%
\alpha \left( u\right) $ as%
\begin{equation*}
\varphi \left( u,v\right) =%
\begin{array}{c}
\left( \cos u\left( 1-\sin \left( \sum_{k=1}^{p}a_{2k}\left( \sinh v\right)
^{k}\right) \right) ,\right.  \\ 
\sin u\left( 1-\sin \left( \sum_{k=1}^{p}a_{2k}\left( \sinh v\right)
^{k}\right) \right) , \\ 
\left. \sin \left( \sum_{k=1}^{p}a_{3k}\left( 1-\cosh v\right) ^{k}\right)
\right) ,%
\end{array}%
\end{equation*}%
where $0\leq v\leq \frac{1}{2}\ \left( Figure\ 3\right) .$%
\begin{figure}[tbh]
\caption{{}A member of spacelike surface family and its common spacelike
asymptotic.}
\end{figure}
\end{example}

\bigskip

Now, we give some special examples. We construct all minimal timelike
surfaces (i.e. helicoid of the 1st, 2nd and 3rd kind and the conjugate
surface of Enneper of the 2nd kind, \cite{woestijne} ) as members of
timelike surface family with a common asymptotic curve.

\begin{example}
\label{example 4.3}(The helicoid of the 1st kind). Let 
\begin{equation*}
\alpha \left( u\right) =\left( \frac{4}{9}\cos 3u,\frac{4}{9}\sin 3u,\frac{5%
}{3}u\right) 
\end{equation*}
be a timelike curve, where $0\leq u\leq 2\pi $. It is easy to show that$\ $%
\begin{equation*}
\left\{ 
\begin{array}{c}
e_{1}\left( u\right) =\left( -\frac{4}{3}\sin 3u,\frac{4}{3}\cos 3u,\frac{5}{%
3}\right) , \\ 
e_{2}\left( u\right) =\left( -4\cos 3u,-4\sin 3u,0\right) , \\ 
e_{3}=\left( \frac{5}{3}\sin 3u,-\frac{5}{3}\cos 3u,-\frac{4}{3}\right) .%
\end{array}%
\right. 
\end{equation*}%
If we choose $x\left( u,v\right) =z\left( u,v\right) =0$ and $y\left(
u,v\right) =\frac{v}{4}\ $ and $v_{0}=0$  then Eqn. \ref{3.3} is satisfied.
Thus, the timelike minimal surface family with common timelike asymptotic is
given by 
\begin{equation*}
\varphi \left( u,v\right) =\left( \left( \frac{4}{9}-v\right) \cos 3u,\left( 
\frac{4}{9}-v\right) \sin 3u,\frac{5}{3}u\right), 
\end{equation*}%
where $-1\leq v\leq 1\ \left( Figure\ 4\right) .\ $This is the
parametrization of the helicoid of the 1st kind.
\begin{figure}[tbh]
\caption{{}A member of timelike minimal surface family and its common
timelike asymptotic (The helicoid of the 1st kind).}
\end{figure}
\end{example}

\begin{example}
\label{example 4.4}(The helicoid of the 2nd kind). Let 
\begin{equation*}
\alpha \left( u\right) =\left( -\frac{5}{9}\cosh 3u,\frac{4}{3}u,-\frac{5}{9}%
\sinh 3u\right) 
\end{equation*}
be a timelike curve, where $-1\leq u\leq 1$. It is easy to show that$\ $%
\begin{equation*}
\left\{ 
\begin{array}{c}
e_{1}\left( u\right) =\left( -\frac{5}{3}\sinh 3u,\frac{4}{3},-\frac{5}{3}%
\cosh 3u\right) , \\ 
e_{2}\left( u\right) =\left( -5\cosh 3u,0,-5\sinh 3u\right) , \\ 
e_{3}=\left( -\frac{4}{3}\sinh 3u,\frac{5}{3},-\frac{4}{3}\cosh 3u\right) .%
\end{array}%
\right. 
\end{equation*}%
If we choose $x\left( u,v\right) =z\left( u,v\right) =0$ and $y\left(
u,v\right) =\frac{v}{5}\ $ and $v_{0}=0$ , then Eqn. \ref{3.3} is satisfied.
Thus, we obtain a member of timelike minimal surface family with common
timelike asymptotic as shown in Figure 5: 
\begin{equation*}
\varphi \left( u,v\right) =\left( \left( -\frac{5}{9}-v\right) \cosh 3u,%
\frac{4}{3}u,\left( -\frac{5}{9}-v\right) \sinh 3u\right), 
\end{equation*}%
where $-1\leq v\leq 1.\ $This is the parametrization of the helicoid of the
2nd kind.
\begin{figure}[tbh]
\caption{A member of timelike minimal surface family and its common timelike
asymptotic (The helicoid of the 2nd kind).}
\end{figure}
\end{example}

\begin{example}
\label{example 4.5}(The helicoid of the 3rd kind). Let 
\begin{equation*}
\alpha \left( u\right) =\left( -\frac{3}{25}\sinh 5u,\frac{4}{5}u,-\frac{3}{%
25}\cosh 5u\right) 
\end{equation*}
be a spacelike curve, where $-1\leq u\leq 1$. It is easy to show that$\ $%
\begin{equation*}
\left\{ 
\begin{array}{c}
e_{1}\left( u\right) =\left( -\frac{3}{5}\cosh 5u,\frac{4}{5},-\frac{3}{5}%
\sinh 5u\right) , \\ 
e_{2}\left( u\right) =\left( -3\sinh 5u,0,-3\cosh 5u\right) , \\ 
e_{3}=\left( -\frac{4}{5}\cosh 5u,-\frac{3}{5},-\frac{4}{5}\sinh 5u\right) .%
\end{array}%
\right. 
\end{equation*}%
If we choose $x\left( u,v\right) =z\left( u,v\right) =0$ and $y\left(
u,v\right) =\frac{v}{3}\ $ and $v_{0}=0$ , then Eqn. \ref{3.3} is satisfied.
Thus, we obtain a member of timelike minimal surface family with common
timelike asymptotic as shown in Figure 6: 
\begin{equation*}
\varphi \left( u,v\right) =\left( \left( -\frac{3}{25}-v\right) \sinh 5u,%
\frac{4}{5}u,\left( -\frac{3}{25}-v\right) \cosh 5u\right) ,
\end{equation*}%
where $-1\leq v\leq 1.\ $This is the parametrization of the helicoid of the
3rd kind.
\begin{figure}[tbh]
\caption{A member of timelike minimal surface family and its common
spacelike asymptotic (The helicoid of the 3rd kind).}
\end{figure}
\end{example}

\begin{example}
\label{example 4.6}(The conjugate surface of Enneper of the 2nd kind). Let%
\begin{equation*}
\alpha \left( u\right) =\left( \frac{u^{2}}{2},-\frac{u3}{6},-\frac{u^{3}}{6}%
+u\right) 
\end{equation*}
be a timelike curve, where $-2\leq u\leq 2$. It is easy to show that$\ $%
\begin{equation*}
\left\{ 
\begin{array}{c}
e_{1}\left( u\right) =\left( u,-\frac{u^{2}}{2},-\frac{u^{2}}{2}+1\right) ,
\\ 
e_{2}\left( u\right) =\left( 1,-u,-u\right) , \\ 
e_{3}=\left( -u,-\frac{u^{2}}{2}-1,-\frac{u^{2}}{2}\right) .%
\end{array}%
\right. 
\end{equation*}%
If we choose $x\left( u,v\right) =z\left( u,v\right) =0$ and $y\left(
u,v\right) =v\ $ and $v_{0}=0$ , then Eqn. \ref{3.3} is satisfied. Thus, we
obtain a member of timelike minimal surface family with common timelike
asymptotic as shown in Figure 7: 
\begin{equation*}
\varphi \left( u,v\right) =\left( \frac{u^{2}}{2}+v,-\frac{u^{3}}{6}-uv,-%
\frac{u^{3}}{6}-uv+u\right) ,
\end{equation*}%
where $-1\leq v\leq 1.$This is the parametrization of the conjugate surface
of Enneper of the 2nd kind.
\begin{figure}[tbh]
\caption{A member of timelike minimal surface family and its common timelike
asymptotic (The conjugate surface of Enneper of the 2nd kind).}
\end{figure}
\end{example}

\end{document}